\documentclass[12pt,reqno]{amsart}
\usepackage{amsfonts}
\usepackage{amsmath}
\usepackage{amsthm}
\usepackage{mathrsfs}
\usepackage{stmaryrd}
\usepackage{graphicx}
\usepackage{epstopdf}
\usepackage{textcomp}

\newtheorem{theorem}{Theorem}

\theoremstyle{definition}

\newtheorem{assumption}{Assumption}

\begin{document}
\title{Central limit theorem under variance uncertainty}

\author{Dmitry B. Rokhlin}

\address{Institute of Mathematics, Mechanics and Computer Sciences,
              Southern Federal University,
Mil'chakova str., 8a, 344090, Rostov-on-Don, Russia}
\email{rokhlin@math.rsu.ru}
\thanks{The research is supported by Southern Federal University, project 213.01-07-2014/07.}

\begin{abstract}
We prove the central limit theorem (CLT) for a sequence of independent zero-mean random variables $\xi_j$, perturbed by predictable multiplicative factors $\lambda_j$ with values in intervals $[\underline\lambda_j,\overline\lambda_j]$. It is assumed that the sequences $\underline\lambda_j$, $\overline\lambda_j$ are bounded and satisfy some stabilization condition. Under the classical Lindeberg condition we show that the CLT limit, corresponding to a ``worst'' sequence $\lambda_j$, is described by the solution $v$ of one-dimensional $G$-heat equation. The main part of the proof follows Peng's approach to the CLT under sublinear expectations, and utilizes H\"{o}lder regularity properties of $v$. Under the lack of such properties, we use the technique of half-relaxed limits from the theory of viscosity solutions.
\end{abstract}

\maketitle

\section{Introduction} \label{sec:1} 
\setcounter{equation}{0}

Consider a sequence of independent one-dimensional random  variables $(\xi_j)_{j=1}^\infty$ with zero means and finite variances $\sigma_j^2=\mathsf E\xi_j^2>0$. Put $s_n^2=\sum_{j=1}^n\sigma_j^2$, $\varepsilon>0$ and assume that the Lindeberg condition
\begin{equation} \label{eq:1.1}
L_n(\varepsilon)=\frac{1}{s_n^2}\sum_{j=1}^n\mathsf E\left(\xi_j^2 I_{\{|\xi_j|>\varepsilon s_n\}}\right)\to 0,\quad n\to\infty
\end{equation}
is satisfied. Then, by the classical central limit theorem (CLT), for any bounded continuous function $f:\mathbb R\mapsto\mathbb R$
we have
\begin{equation} \label{eq:1.2}
\lim_{n\to\infty}\mathsf E f\left(\frac{1}{s_n}\sum_{j=1}^n \xi_j\right)=\mathsf E f(\zeta),
\end{equation}
where $\zeta$ has the standard normal law.

In this paper we assume that the variance of $\xi_j$ is not known exactly and may belong to an interval.
Our goal is to obtain the ``least upper bound'' $\mathscr L$ for the quantity (\ref{eq:1.2}) under such model uncertainty. The result, as well as its proof, are similar to those obtained by Peng \cite{Peng07,Peng08} and the followers \cite{LiShi10,ZhaChe14,HuZhou15} under the nonlinear expectations theory paradigm. It appears that $\mathscr L$ can be described in terms of the solution $v$ of a nonlinear parabolic equation, called $G$-heat equation. One of the objectives of the present paper is to show that this  description also comes from a classical problem statement, and need not be linked to the nonlinear expectations theory.

To give a precise problem formulation, consider a filtered probability space
\[(\Omega,\mathscr F,\mathbb P, (\mathscr F_j)_{j=0}^\infty)\]
and an adapted sequence $(\xi_j)_{j=1}^\infty$ of random variables such that $\mathsf E\xi_j=0$, $\mathsf E\xi_j^2=\sigma_j^2\in (0,\infty)$  and $\xi_j$ is independent from $\mathscr F_{j-1}$. Let $(\lambda_j)_{j=0}^\infty$ be an adapted sequence, whose elements $\lambda_j$ take values in deterministic intervals $[\underline\lambda_j,\overline\lambda_j]$, $0\le\underline\lambda_j\le\overline \lambda_j$. Considering the sequence $\eta_j=\lambda_{j-1}\xi_j$, one can regard the multipliers $\lambda_{j-1}$ as a ``predictable perturbation'' of the original sequence $\xi_j$.
The intervals $[\underline\lambda_{j-1}\sigma_j,\overline \lambda_{j-1}\sigma_j]$ indicate possible standard deviations of $\eta_j$.

\begin{assumption} \label{as:1}
The Lindeberg condition (\ref{eq:1.1}) is satisfied.
\end{assumption}
\begin{assumption} \label{as:2}
The sequence $\overline \lambda_j$ is bounded by a constant $\Lambda$.
\end{assumption}
\begin{assumption} \label{as:3}
The sequences $\underline\lambda_j$, $\overline\lambda_j$ satisfy the following \emph{stabilization condition}:
\begin{equation} \label{eq:1.3}
M_n=\sum_{j=0}^{n-1}\frac{\sigma_{j+1}^2}{s_n^2}\left(|\overline\lambda_j^2-\overline\lambda^2|+
|\underline\lambda_j^2-\underline\lambda^2|\right)\to 0,\quad n\to\infty
\end{equation}
for some $\overline\lambda\ge\underline\lambda\ge 0$.
\end{assumption}

Put $B_j=[\underline\lambda_j^2,\overline\lambda_j^2]$, $B=[\underline\lambda^2,\overline\lambda^2]$ and denote by
$$ d_H(B_j,B)=\max\{|\overline\lambda_j^2-\overline\lambda^2|,|\underline\lambda_j^2-\underline\lambda^2|\}$$
the Hausdorff distance between these intervals (see, e.g., \cite[Chapter 2]{AleHer83}). Condition (\ref{eq:1.3}) is equivalent to the following one:
\begin{equation} \label{eq:1.3A}
\sum_{j=0}^{n-1}\frac{\sigma_{j+1}^2}{s_n^2}d_H(B_j,B)\to 0,\quad n\to\infty.
\end{equation}
In the summability theory the transformation
\[ t_n=\frac{p_1 a_1+\dots+p_n a_n}{p_1+\dots+p_n},\quad p_n>0 \]
of a sequence $(a_i)_{i=1}^\infty$ is called a Riesz mean (see \cite[Section 1.4]{Pet66}, \cite[Section 3.2]{Boos00}). By Assumption \ref{as:3}, the sequence $d_H(B_j,B)$ is summable to $0$ by the Riesz method, determined by the sequence $p_i=\sigma_i^2$. Furthermore, the Lindeberg condition implies the Feller condition
\begin{equation} \label{eq:1.3B}
\lim_{n\to\infty}\max_{1\le j\le n}\frac{\sigma_j}{s_n}=0
\end{equation}
(see, e.g., \cite{Bau96}, Chapter 6, \S 28). In particular, $s_n\to\infty$. Hence, the Riesz summation method, defined above, is regular (see \cite[Theorem 1.4.4]{Pet66}), and if $d_H(B_j,B)\to 0$, then the Assumption \ref{as:3} is satisfied. We also mention a necessary and sufficient condition for (\ref{eq:1.3A}) to hold true, given in \cite{Boos00} (Lemma 3.2.14). This result is applicable since $\sigma_n^2/s_n^2\to 0$ by (\ref{eq:1.3B}).

Note that from the identity
\[ \sum_{j=0}^{n-1}\frac{\sigma_{j+1}^2}{s_n^2}=1 \]
it easily follows that
\[\lim_{n\to\infty}\sum_{j=0}^{n-1}\frac{\sigma_{j+1}^2}{s_n^2}\overline\lambda_j^2=\overline\lambda^2, \quad
   \lim_{n\to\infty}\sum_{j=0}^{n-1}\frac{\sigma_{j+1}^2}{s_n^2}\underline\lambda_j^2=\underline\lambda^2. \]

Denote by $\mathfrak A^n$ the set of adapted sequences $\lambda_0^n=(\lambda_j)_{j=0}^n$ with values in $[\underline\lambda_j,\overline\lambda_j]$. Our goal is to describe the quantity
\begin{equation} \label{eq:1.5}
\mathscr L=\lim_{n\to\infty}\sup_{\lambda_0^{n-1}\in\mathfrak A^{n-1}}\mathsf E f\left(\frac{1}{s_n}\sum_{j=0}^{n-1} \lambda_j\xi_{j+1}\right),
\end{equation}
which can be loosely characterized as the least upper bound of (\ref{eq:1.2}) under variance uncertainty.

The main role in this description is played by the solution of the nonlinear parabolic equation
\begin{equation} \label{eq:1.6}
v_t+\frac{1}{2}\sup_{\lambda\in[\underline\lambda,\overline\lambda]}\left(\lambda^2 v_{xx}\right)=
v_t+\frac{1}{2}\left(\overline\lambda^2 v_{xx}^+-\underline\lambda^2v_{xx}^- \right)=0, \quad (t,x)\in [0,1)\times\mathbb R,
\end{equation}
satisfying the boundary condition
\begin{equation} \label{eq:1.7}
v(1,x)=f(x),\quad x\in\mathbb R.
\end{equation}

In the context of the CLT under sublinear expectations, equation (\ref{eq:1.6}) appeared in \cite{Peng07}. It was called $G$-heat equation in \cite{Peng07b}. As is mentioned in \cite{CafSte08}, such equation arises in various applications in control theory, mechanics, combustion, biology, and finance. It is known also as a Barenblatt equation: see, e.g., \cite{KamPelVaz91}.

One can obtain (\ref{eq:1.6}) by considering $\lambda_j$ as a control sequence, writing down dynamic programming equations for discrete time finite horizon optimization problems
\[ \sup_{\lambda_0^{n-1}\in\mathfrak A^{n-1}}\mathsf E f\left(\frac{1}{s_n}\sum_{j=0}^{n-1} \lambda_j\xi_{j+1}\right),\]
and passing to the limit as $n\to\infty$. This approach was proposed in \cite{Rok15} in the case of identically distributed (multidimensional) random variables $\xi_j$. However, in the present context, it seems that this method requires hypotheses, which are stronger than the Lindeberg condition. Thus, we follow Peng's approach, which takes equation (\ref{eq:1.6}) as a starting point, and utilizes a deep result on the existence of its solution in an appropriate H\"older class.

Put $Q=[0,1]\times\mathbb R$,
\[ \|h\|_{0;\mathbb R}=\sup_{x\in\mathbb R} |h(x)|,\quad \|g\|_{0;Q}=\sup_{(t,x)\in Q} |g(t,x)|, \]
\[ [h]_{\alpha;\mathbb R}=\sup_{\substack{x_i\in\mathbb R,\\ x_1\neq x_2}} \frac{|h(x_1)-h(x_2)|}{|x_1-x_2|^\alpha},\quad \alpha\in (0,1], \]
\[ [g]_{\alpha;Q}=\sup_{\substack{(t_i,x_i)\in Q,\\ (t_1,x_1)\neq (t_2,x_2)}} \frac{|g(t_1,x_1)-g(t_2,x_2)|}{(|t_1-t_2|^{1/2}+|x_1-x_2|)^\alpha},\quad
\alpha\in (0,1], \]
and consider the H\"older spaces $C^{2+\alpha}(\mathbb R)$, $C^{1+\alpha/2,2+\alpha}(Q)$ with the norms
\[ \|h\|_{C^{2+\alpha}(\mathbb R)}=\|h\|_{0;\mathbb R}+\|h_x\|_{0;\mathbb R}+\|h_{xx}\|_{0;\mathbb R}+[h_{xx}]_{\alpha;\mathbb R},\]
\[ \|g\|_{C^{1+\alpha/2,2+\alpha}(Q)}=\|g\|_{0;Q}+\|g_x\|_{0;Q}+\|g_t\|_{0;Q}+\|g_{xx}\|_{0;Q}+[g_t]_{\alpha;Q}+ [g_{xx}]_{\alpha;Q}.\]

Under the assumptions $f\in C^{2+\alpha}(\mathbb R)$, $\alpha\in (0,1]$; $\underline\lambda>0$ the existence of a classical solution $v\in C^{1+\alpha'/2,2+\alpha'}(Q)$ (with some of $\alpha'\in (0,1]$) of (\ref{eq:1.6}), (\ref{eq:1.7}) was proved by Krylov: see \cite{Kry82} (Theorem 1.1 or Theorem 5.3).

If $\underline\lambda=0$ then only the existence of a viscosity solution is guaranteed. Let us recall this result along with related definitions. Put $Q^\circ=[0,1)\times\mathbb R$ and assume that $f$ is a bounded continuous function: $f\in C_b(\mathbb R)$. A bounded upper semicontinuous (usc) function $v:Q\mapsto\mathbb R$ is called a \emph{viscosity subsolution} of (\ref{eq:1.6}), (\ref{eq:1.7}) if
\begin{equation} \label{eq:1.8}
v(1,x)\le f(x),\quad x\in\mathbb R,
\end{equation}
and for any $(\overline t, \overline x)\in Q^\circ$ and any test function $\varphi\in C^2(\mathbb R^2)$ such that $(\overline t, \overline x)\in Q^\circ$ is a strict local maximum point of $v-\varphi$ on $Q^\circ$, the inequality
\begin{equation} \label{eq:1.9}
-\varphi_t(\overline t, \overline x)-\frac{1}{2}\sup_{\lambda\in[\underline\lambda,\overline\lambda]}\left(\lambda^2 \varphi_{xx}(\overline t, \overline x)\right)\le 0
\end{equation}
holds true. To define a \emph{viscosity supersolution}, one should consider a bounded lower semicontinuous (lsc) function $v$, a strict local minimum point of $v-\varphi$, and reverse the inequalities (\ref{eq:1.8}), (\ref{eq:1.9}).

We will use the following comparison result. Consider a viscosity subsolution $u$ and a viscosity supersolution $w$ of (\ref{eq:1.6}), (\ref{eq:1.7}). Since we require (\ref{eq:1.6}) to be satisfied in the viscosity sense at the lower boundary of $Q$, by the accessibility theorem of \cite{CheGigGot91}, we have
\[ u(0,x)=\limsup_{\substack{(t,y)\in (0,1)\times\mathbb R,\\ t\to 0, y\to x}} u(t,y);\quad
   w(0,x)=\liminf_{\substack{(t,y)\in (0,1)\times\mathbb R,\\ t\to 0, y\to x}} w(t,y)
\]
and by the comparison result of \cite{DieFriObe14} (Theorem 1) it follows that $u\le w$ on $Q$.

A bounded \emph{continuous} function $v:Q\mapsto\mathbb R$ is called a \emph{viscosity solution} of (\ref{eq:1.6}), (\ref{eq:1.7}), if it is viscosity sub- and supersolution. The existence of a continuous viscosity solution of (\ref{eq:1.6}), (\ref{eq:1.7}) for $f\in C_b(\mathbb R)$ is well known from the theory of optimal control. The stochastic control representation of such solution can be found in \cite{YongZhou99} (Chap. 4, Theorem 5.2).

\begin{theorem} \label{th:1}
Let $f$ be a bounded continuous function, and let $v$ be the continuous viscosity solution of (\ref{eq:1.6}), (\ref{eq:1.7}). Then, under Assumptions \ref{as:1}--\ref{as:3}, we have $\mathscr L=v(0,0)$.
\end{theorem}

It is interesting to compare Theorem \ref{th:1} with related results obtained in the framework of sublinear expectations theory. Besides the original result of Peng \cite{Peng07,Peng08}, which is discussed in \cite{Rok15}, we mention the papers \cite{LiShi10,ZhaChe14,HuZhou15}, where the random variables were not assumed to be identically distributed. We will discuss only the result of \cite{ZhaChe14}, which extends \cite{LiShi10}. The result of \cite{HuZhou15} concerns the multidimensional case.

Let us briefly describe the construction of a sublinear expectation space $(\Omega,\mathcal H,\widehat{\mathbb E})$, which allows to rewrite the expression (\ref{eq:1.5}) in terms of a sublinear expectation. This construction is, in fact, the same as in \cite[Section 4]{Rok15}, where some more details are given. Consider the space of sequences $\Omega=\{(y_i)_{i=1}^\infty:y_i\in\mathbb R\}$, and introduce the space of random variables $\mathcal H$ as follows: $\mathcal H=\cup_{n=1}^\infty \mathcal H_n$, where
$\mathcal H_n$ is some linear  space (we do not go into details) of  functions $Y=\psi(y_1,\dots,y_n)$ of $n$ variables. Define the sublinear expectation by the formula
\begin{equation} \label{eq:1.11}
 \widehat{\mathbb E}Y=\sup_{\lambda_0^{n-1}\in\mathfrak A_0^{n-1}}\mathsf E\psi(\lambda_0\xi_1/\sigma_1,\dots,\lambda_{n-1}\xi_n/\sigma_n).
\end{equation}
Let $Y_i$ be the projection mappings: $Y_i(y)=y_i$. One can show that $Y_n$ is independent from $(Y_1,\dots,Y_{n-1})$ in the sense of sublinear expectations theory (see \cite{Peng08c}, Definition 3.10).
By (\ref{eq:1.11}) we get the following representation for $\mathscr L$:
$$ \mathscr L=\lim_{n\to\infty}\sup_{\lambda_0^{n-1}\in\mathfrak A^{n-1}}\mathsf E f\left(\frac{1}{s_n}\sum_{j=1}^{n} \lambda_{j-1}\xi_{j}\right)=\lim_{n\to\infty}\widehat{\mathbb E} f\left(\sum_{j=1}^n\frac{\sigma_j}{s_n} Y_j\right).$$

Let us apply Theorem 3.1 of \cite{ZhaChe14} to the sequence $Y_i$. We have
\begin{align*}
 \widehat{\mathbb E}(\pm Y_i) &= \sup_{\lambda_{i-1}\in [\underline\lambda_{i-1},\overline\lambda_{i-1}]}\mathsf E(\pm\lambda_{i-1}\xi_i/\sigma_i)=0,\\
 \widehat{\mathbb E} Y_i^2 &=\sup_{\lambda_{i-1}\in [\underline\lambda_{i-1},\overline\lambda_{i-1}]}\mathsf E(\lambda_{i-1}\xi_i/\sigma_i)^2 =\overline\lambda_{i-1}^2, \\
-\widehat{\mathbb E}(-Y_i^2) &=-\sup_{\lambda_{i-1}\in [\underline\lambda_{i-1},\overline\lambda_{i-1}]}\mathsf E\left(-(\lambda_{i-1}\xi_i/\sigma_i)^2\right)=\underline\lambda_{i-1}^2.
\end{align*}
Besides a condition, identical to Assumption \ref{as:3}, in \cite{ZhaChe14} it is assumed that
\begin{equation} \label{eq:1.12}
 \widehat{\mathbb E}|Y_i|^{2+\delta}=\sup_{\lambda_{i-1}\in [\underline\lambda_{i-1},\overline\lambda_{i-1}]}\mathsf E\left|\lambda_{i-1}\xi_i/\sigma_i\right|^{2+\delta}=\overline\lambda_{i-1}^{2+\delta}\mathsf E|\xi_i/\sigma_i|^{2+\delta}\le M,
\end{equation}
\begin{equation} \label{eq:1.13}
 \lim_{n\to\infty} \sum_{j=1}^n\left(\frac{\sigma_j}{s_n}\right)^{2+\delta}=0
\end{equation}
for some $M>0$, $\delta>0$. Note that (\ref{eq:1.12}) was used in \cite{ZhaChe14} in this form, although it was not clearly formulated (see condition (3) of Theorem 3.1 in \cite{ZhaChe14}). The result of \cite{ZhaChe14} tells us that $\mathscr L=\widehat{\mathbb E} f(Z)$, where $Z$ is a $G$-normal random variable with
$$ G(s)=\frac{1}{2}(s^+\overline\lambda^2-s^-\underline\lambda^2).$$
By the characterization of the $G$-normal distribution (see, e.g., \cite{Peng08c}, Example 1.13) this is equivalent to the assertion of Theorem \ref{th:1}.

Thus, under the assumptions (\ref{eq:1.12}), (\ref{eq:1.13}) (instead of Assumptions \ref{as:1}, \ref{as:2}), Theorem \ref{th:1} follows from the result of \cite{ZhaChe14}. It is easy to see that (\ref{eq:1.12}) implies Assumption \ref{as:2}:
\[ \overline\lambda_{i-1}=\left(\mathsf E(\overline\lambda_{i-1}^2\xi_i^2/\sigma_i^2)\right)^{1/2}
\le\left(\mathsf E(\overline\lambda_{i-1}^{2+\delta}|\xi_i/\sigma_i|^{2+\delta})\right)^{1/(2+\delta)} \]
and (\ref{eq:1.12}), (\ref{eq:1.13}) imply that
\[\frac{1}{s_n^2}\sum_{j=1}^n\mathsf E\left(\overline\lambda_{j-1}^{2+\delta}\xi_j^2 I_{\{|\xi_j|>\varepsilon s_n\}}\right)\le
\frac{1}{\varepsilon^\delta s_n^{2+\delta}}\sum_{j=1}^n\overline\lambda_{j-1}^{2+\delta}\mathsf E\left|\xi_j\right|^{2+\delta}
\le\frac{M}{\varepsilon^\delta s_n^{2+\delta}}\sum_{j=1}^n\sigma_j^{2+\delta}\to 0,\quad n\to\infty. \]
The last condition is slightly weaker than Assumption \ref{as:1}, and coincides with the latter if $\liminf_{j\to\infty} \overline\lambda_j>0$.

Note that if there is no model uncertainty: $\underline\lambda_j=\overline\lambda_j=1$, then Theorem \ref{th:1} reduces to the classical CLT, mentioned at the beginning of the present paper. This is not the case with the result of \cite{ZhaChe14}, since in this case the conditions (\ref{eq:1.12}), (\ref{eq:1.13}) are stronger then the Lindeberg condition. We also mention that \cite{ZhaChe14} deals only with classical solutions of the $G$-heat equation, so the case $\underline\lambda=0$ is, in fact, not considered there. However, the sublinear expectations theory is able to handle the degenerate case via perturbation methods, see \cite{Peng08} (the proof of Theorem 5.1), \cite{HuZhou15} (the proof of Theorem 3.1).

\section{Proof of Theorem \ref{th:1}}

(i) We first consider the case $f\in C^{2+\alpha}(\mathbb R)$, $\alpha>0$ and $\underline\lambda>0$. Put
\[ X_{j+1}=X_j+\frac{\lambda_j}{s_n}\xi_{j+1},\quad j=0,\dots,n-1,\quad X_0=0;\quad t_j=\sum_{k=0}^j \frac{\sigma_k^2}{s_n^2}.\]
Since the solution $v$ of (\ref{eq:1.6}), (\ref{eq:1.7}) belongs to $v\in C^{1+\alpha'/2,2+\alpha'}(Q)$, we can apply Taylor's formula:
\begin{align*}
& v(1,X_n)-v(0,0)=\sum_{j=0}^{n-1}\left(v(t_{j+1},X_{j+1})-v(t_j,X_{j+1})+
v(t_j,X_{j+1})-v(t_j,X_j)\right)\\
=& \sum_{j=0}^{n-1}\left(v_t(\widehat t_j,X_{j+1})(t_{j+1}-t_j) +v_x(t_j,X_j)(X_{j+1}-X_j)+\frac{1}{2}v_{xx}(t_j,\widehat X_j)(X_{j+1}-X_j)^2\right),
\end{align*}
where $\widehat t_j=t_j+\beta (t_{j+1}-t_j)$, $\widehat X_j=X_j+\gamma (X_{j+1}-X_j)$, $\beta,\gamma\in [0,1]$. By the independence of $X_j$ and $\xi_{j+1}$ we conclude that  $\mathsf E(v_x(t_j,X_j)(X_{j+1}-X_j))=0$. Thus,
\[\mathsf E v(1,X_n)-v(0,0)=\mathsf E\sum_{j=0}^{n-1}\frac{\sigma_{j+1}^2}{s_n^2}\left(v_t(\widehat t_j,X_{j+1}) + \frac{\lambda_j^2}{2}\frac{\xi_{j+1}^2}{\sigma_{j+1}^2} v_{xx}(t_j,\widehat X_j)\right)=J_n+I_n,\]
\[ J_n=\mathsf E\sum_{j=0}^{n-1}\frac{\sigma_{j+1}^2}{s_n^2}\left(v_t(t_j,X_j) + \frac{\lambda_j^2}{2}\frac{\xi_{j+1}^2}{\sigma_{j+1}^2} v_{xx}(t_j,X_j)\right)=\mathsf E\sum_{j=0}^{n-1}\frac{\sigma_{j+1}^2}{s_n^2}\left(v_t + \frac{\lambda_j^2}{2}v_{xx}\right)(t_j,X_j), \]
\[ I_n=\mathsf E\sum_{j=0}^{n-1}\frac{\sigma_{j+1}^2}{s_n^2}\left(v_t(\widehat t_j,X_{j+1})-v_t(t_j,X_j) + \frac{\lambda_j^2}{2}\frac{\xi_{j+1}^2}{\sigma_{j+1}^2} (v_{xx}(t_j,\widehat X_j)-v_{xx}(t_j,X_j))\right). \]

We can rewrite $J_n$ as $J_n^1+J_n^2$, where
\begin{align*}
 J_n^1 &=\mathsf E\sum_{j=0}^{n-1}\frac{\sigma_{j+1}^2}{s_n^2}\left(v_t + \frac{1}{2}\left( \overline\lambda^2 v_{xx}^+ -\underline\lambda^2 v_{xx}^-\right)\right)(t_j,X_j),\\
 J_n^2 &=\frac{1}{2}\mathsf E\sum_{j=0}^{n-1}\frac{\sigma_{j+1}^2}{s_n^2}\left((\lambda_j^2-\overline\lambda^2)v_{xx}^+ +(\underline\lambda^2-\lambda_j^2) v_{xx}^- \right)(t_j,X_j).
\end{align*}

From the definition of $v$ we see that $J_n^1=0$. Furthermore, from the stabilization condition (\ref{eq:1.3}) it follows that
\[ \sup_{\lambda_0^{n-1}\in\mathfrak A^{n-1}}J_n^2\le\frac{1}{2}\mathsf E \sum_{j=0}^{n-1}\frac{\sigma_{j+1}^2}{s_n^2}\left((\overline \lambda_j^2-\overline\lambda^2)v_{xx}^+ +(\underline\lambda^2-\underline\lambda_j^2) v_{xx}^-\right)(t_j,X_j)\le C M_n\to 0, \]
$n\to\infty$, since the second derivative of $v$ is uniformly bounded. On the other hand, choosing a sequence
\[ \lambda_j=\overline\lambda_j I_{\{v_{xx}(t_j,X_j)>0\}}-\underline\lambda_j I_{\{v_{xx}(t_j,X_j)\le 0\}},\quad j\ge 1,  \]
with an arbitrary $\lambda_0$, we get an opposite inequality
\[ \sup_{\lambda_0^{n-1}\in\mathfrak A^{n-1}}J_n^2\ge \frac{1}{2}\mathsf E \sum_{j=0}^{n-1}\frac{\sigma_{j+1}^2}{s_n^2}\left((\overline \lambda_j^2-\overline\lambda^2)v_{xx}^+ +(\underline\lambda^2-\underline\lambda_j^2) v_{xx}^- \right)(t_j,X_j)\ge -C M_n\to 0. \]
Combining all these results, we conclude that
\begin{equation} \label{eq:2.1}
\lim_{n\to\infty}\sup_{\lambda_0^{n-1}\in\mathfrak A^{n-1}} J_n=0.
\end{equation}

Now consider $I_n=I_n^1+I_n^2+I_n^3$:
\[ I_n^1=\mathsf E\sum_{j=0}^{n-1}\frac{\sigma_{j+1}^2}{s_n^2}\left(v_t(\widehat t_j,X_{j+1})-v_t(t_j,X_j)\right), \]
\[ I_n^2=\mathsf E\sum_{j=0}^{n-1}\frac{\xi_{j+1}^2}{s_n^2}\frac{\lambda_j^2}{2}\left( v_{xx}(t_j,\widehat X_j)-v_{xx}(t_j,X_j)\right)I_{\{|\xi_{j+1}|>\varepsilon s_n\}}, \]
\[ I_n^3=\mathsf E\sum_{j=0}^{n-1}\frac{\xi_{j+1}^2}{s_n^2}\frac{\lambda_j^2}{2}\left( v_{xx}(t_j,\widehat X_j)-v_{xx}(t_j,X_j)\right)I_{\{|\xi_{j+1}|\le\varepsilon s_n\}}. \]

By the H\"{o}lder continuity of $v_t$ we have
\begin{align*}
 |I_n^1| \le & C\mathsf E\sum_{j=0}^{n-1}\frac{\sigma_{j+1}^2}{s_n^2}\left(|\widehat t_j-t_j|^{\alpha'/2}+|X_{j+1}-X_j|^{\alpha'}\right)\\
 \le & C\mathsf E\sum_{j=0}^{n-1}\frac{\sigma_{j+1}^2}{s_n^2}\left(\left(\frac{\sigma_{j+1}}{s_n}\right)^{\alpha'}+ \left(\frac{\lambda_j|\xi_{j+1}|}{s_n}\right)^{\alpha'} \right).
\end{align*}
Using the inequality $\mathsf E|\xi_{j+1}|^{\alpha'}\le(\mathsf E\xi_{j+1}^2)^{{\alpha'}/2}=\sigma_{j+1}^{\alpha'}$, and the independence of $\lambda_j$ and $\xi_{j+1}$, we obtain the estimate
\[ |I_n^1| \le C \sum_{j=0}^{n-1}\frac{\sigma_{j+1}^2}{s_n^2}(1+\overline\lambda_j^{\alpha'})
\left(\frac{\sigma_{j+1}}{s_n}\right)^{\alpha'}\le C \left(\max_{1\le j\le n}\frac{\sigma_j}{s_n}\right)^{\alpha'} \left(1+
\Lambda^{\alpha'}\right). \]
From (\ref{eq:1.3B}) it follows that $I_n^1\to 0$.

Furthermore, since the sequence $\lambda_j$ is bounded and the second derivative of $v$ is uniformly bounded, by the Lindeberg condition we get
\[ |I_n^2|\le C L_n(\varepsilon)\to 0,\quad n\to\infty. \]
The last term $I_n^3$ is estimated with the use of the H\"{o}lder continuity property of $v_{xx}$:
\[ |I_n^3|\le C \mathsf E\sum_{j=0}^{n-1}\frac{\xi_{j+1}^2}{s_n^2}\frac{\lambda_j^2}{2} \left|\frac{\lambda_j |\xi_{j+1}|}{s_n}\right|^{\alpha'} I_{\{|\xi_{j+1}|\le\varepsilon s_n\}}\le C\Lambda^{2+\alpha'}\sum_{j=0}^{n-1}\frac{\sigma_{j+1}^2}{s_n^2}\varepsilon^{\alpha'}=C\Lambda^{2+\alpha'}\varepsilon^{\alpha'}. \]

Therefore,
\begin{equation} \label{eq:2.2}
 \lim_{n\to\infty}\sup_{\lambda_0^{n-1}\in\mathfrak A^{n-1}} |I_n|=0.
\end{equation}

From (\ref{eq:2.1}) and (\ref{eq:2.2}) it follows that
$$\mathscr L=\lim_{n\to\infty}\sup_{\lambda_0^{n-1}\in\mathfrak A^{n-1}}\mathsf E f(X_n)=\lim_{n\to\infty}\sup_{\lambda_0^{n-1}\in\mathfrak A^{n-1}}\mathsf E v(1,X_n)=v(0,0).$$
So, we have proved the theorem in the case $f\in C^{2+\alpha}$, $\underline\lambda>0$.

(ii) Now assume that $\underline\lambda=0$. Put
$$ X_n^\varepsilon=\frac{1}{s_n}\sum_{j=0}^{n-1} (\lambda_j^2+\varepsilon^2)^{1/2}\xi_{j+1}, \quad\mathscr L^\varepsilon=\lim_{n\to\infty}\sup_{\lambda_0^{n-1}\in\mathfrak A^{n-1}}\mathsf E f\left(X_n^\varepsilon\right).$$
The intervals $[\underline\mu_j,\overline\mu_j]=[(\underline\lambda_j^2+\varepsilon^2)^{1/2},(\overline\lambda_j^2+\varepsilon^2)^{1/2}]$ stabilize to $[\varepsilon,(\overline\lambda^2+\varepsilon^2)^{1/2}]$ in the sense of Assumption \ref{as:3}:
$$\sum_{j=0}^{n-1}\frac{\sigma_{j+1}^2}{s_n^2}\left(|\overline\mu_j^2-(\overline\lambda^2+\varepsilon^2)|+
|\underline\mu_j^2-\varepsilon^2|\right)\to 0,\quad n\to\infty.$$

By part (i) of the proof, we infer that $\mathscr L^\varepsilon=v^\varepsilon(0,0)$, where $v^\varepsilon$ satisfies
\begin{equation} \label{eq:2.3}
 v^\varepsilon_t+\frac{1}{2}\left((\overline\lambda^2+\varepsilon^2) (v^\varepsilon_{xx})^+-\varepsilon^2 (v^\varepsilon_{xx})^- \right)=0,\quad x\in Q^\circ; \quad v^\varepsilon(1,x)=f(x),\quad x\in\mathbb R
\end{equation}
in the classical sense. Let $v$ be the continuous viscosity solution of the limiting problem
\begin{equation} \label{eq:2.4}
 v_t+\frac{1}{2}\overline\lambda^2 v_{xx}^+ =0,\quad x\in Q^\circ;\quad v(1,x)=f(x),\quad x\in\mathbb R.
\end{equation}
The desired result is a consequence of the relations
\begin{equation} \label{eq:2.5}
 \mathscr L:=\lim_{n\to\infty}\sup_{\lambda_0^{n-1}\in\mathfrak A^{n-1}}\mathsf E f(X_n) =\lim_{\varepsilon\to 0}\mathscr L^\varepsilon,\quad v(0,0)=\lim_{\varepsilon\to 0} v^\varepsilon(0,0),
\end{equation}
which we are going to prove.

Since we still assume that $f\in C^{2+\alpha}(\mathbb R)$, this function is uniformly Lipschitz continuous. Put $\psi_\varepsilon(\lambda)=(\lambda^2+\varepsilon^2)^{1/2}-\lambda$. We have
\begin{align*}
 & |\mathsf E f(X_n^\varepsilon)-\mathsf E f(X_n)|\le C\mathsf E |X_n^\varepsilon-X_n|\le \frac{C}{s_n}\left(\mathsf E \left(\sum_{j=0}^{n-1}\psi_\varepsilon(\lambda_j)\xi_{j+1}\right)^2\right)^{1/2}\\
 &= C\left(\mathsf E\sum_{j=0}^{n-1}\frac{\sigma_{j+1}^2}{s_n^2}\psi_\varepsilon^2 (\lambda_j)\right)^{1/2}\le C\varepsilon,
\end{align*}
since $\sup_{\lambda\ge 0}\psi_\varepsilon(\lambda)=\varepsilon$. Thus,
\begin{equation} \label{eq:2.6}
\mathscr L^\varepsilon-C\varepsilon\le\liminf_{n\to\infty}\sup_{\lambda_0^{n-1}\in\mathfrak A^{n-1}}\mathsf E f(X_n)
\le \limsup_{n\to\infty}\sup_{\lambda_0^{n-1}\in\mathfrak A^{n-1}}\mathsf E f(X_n)\le \mathscr L^\varepsilon+C\varepsilon,
\end{equation}
\[ \limsup_{\varepsilon\to 0}\mathscr L^\varepsilon\le\liminf_{n\to\infty}\sup_{\lambda_0^{n-1}\in\mathfrak A^{n-1}}\mathsf E f(X_n)
\le \limsup_{n\to\infty}\sup_{\lambda_0^{n-1}\in\mathfrak A^{n-1}}\mathsf E f(X_n)\le \liminf_{\varepsilon\to 0}\mathscr L^\varepsilon. \]
These estimates imply the first equality in (\ref{eq:2.5}).

Furthermore, define the \emph{half-relaxed} (or \emph{weak}) limits of $v^\varepsilon$ by
\[ \underline v(t,x)=\liminf_{\substack{(s,y)\to(t,x),\\ \varepsilon\to 0}} v^\varepsilon(s,y),\quad
\overline v(t,x)=\limsup_{\substack{(s,y)\to(t,x),\\ \varepsilon\to 0}} v^\varepsilon(s,y),\quad (t,x)\in Q.
\]
The function $\overline v$ (resp., $\underline v$) is usc (resp., lsc): see \cite{BarCap97} (Chap.\,5, Lemma 1.5).

Take $\varphi\in C^2(\mathbb R^2)$ and assume that $\overline z=(\overline t,\overline x)\in Q$ is a strict local maximum point of $\overline v-\varphi$ on $Q$. Then there exist sequences $\varepsilon_k\to 0$, $z_k=(t_k,x_k)\in Q$ such that $z_k\to \overline z$, $v^{\varepsilon_k}(z_k)\to\overline v(\overline z)$, and $z_k$ is a local maximum point of $v^{\varepsilon_k}-\varphi$ on $Q$: see \cite{BarCap97} (Chap.\,5, Lemma 1.6).

If $t\in [0,1)$, then $t_k\in [0,1)$ for sufficiently large $k$ and
\[ -\varphi_t(z_k)-\sup_{\lambda\in [\varepsilon_k,\overline\lambda+\varepsilon_k]} (\lambda^2\varphi_{xx}(z_k))\le 0, \]
since $v^{\varepsilon_k}$ is a viscosity solution of (\ref{eq:2.3}). Passing to the limit as $\varepsilon_k\to 0$, we get the inequality
\begin{equation} \label{eq:2.7}
-\varphi_t(\overline z)-\sup_{\lambda\in [0,\overline\lambda]} (\lambda^2\varphi_{xx}(\overline z))\le 0,
\end{equation}
which means that $\overline v$ is a viscosity subsolution of (\ref{eq:2.4}) on $Q^\circ$.

Let $\overline t=1$. If there are infinitely many $t_k<1$, then we again obtain (\ref{eq:2.7}) as above. Moreover, we can change the test function $\varphi$ to $\widehat\varphi=\varphi+c(1-t)$, $c>0$ since $(1,\overline x)$ is still a strict local maximum point of $\overline v-\widehat\varphi$. Substituting  $\widehat\varphi$ in (\ref{eq:2.7}), we get a contradiction:
\[ c-\varphi_t(\overline z)-\sup_{\lambda\in [0,\overline\lambda]} (\lambda^2\varphi_{xx}(\overline z))\le 0,\quad \text{for any } c>0. \]
Thus, for sufficiently large $k$, we have $v^{\varepsilon_k}(z_k)=f(x_k)$ and $\overline v(\overline z)=\lim_{k\to\infty} f(x_k)=f(\overline x)$.

We have proved that $\overline v$ is a viscosity subsolution of (\ref{eq:2.4}). Similarly, one can prove that  $\underline v$ is a viscosity supersolution of (\ref{eq:2.4}). By the comparison result of \cite{DieFriObe14}, mentioned in Section \ref{sec:1}, we have $\overline v\le\underline v$ on $Q$. The converse inequality $\overline v\ge\underline v$ is immediate from the definition. We infer that $v=\overline v=\underline v$ is a continuous viscosity solution of (\ref{eq:2.4}), and the second equality in (\ref{eq:2.5}) holds true:
\[ v(0,0)\le\lim\inf_{\varepsilon\to 0} v^\varepsilon(0,0)\le\limsup_{\varepsilon\to 0} v^\varepsilon(0,0)\le v(0,0). \]
This completes the proof of Theorem \ref{th:1} in the case $\underline\lambda=0$.

(iii) It remains to consider the case $f\in C_b(\mathbb R)$. It is not difficult to show that there exists a function $f^\varepsilon\in C^\infty(\mathbb R)$ such that $|f(x)-f^\varepsilon(x)|\le\varepsilon$: see, e.g., \cite{Pur67}. Furthermore, consider a function $\chi\in C^\infty$,
\[ \chi(x)=1,\quad |x|\le 1;\quad \chi(x)=0,\quad |x|\ge 2, \]
and put $g^\varepsilon(x)=\chi(\varepsilon^{1/2} x)f^\varepsilon(x)$. We have
\begin{align*}
& |\mathsf E f(X_n)-\mathsf E g^\varepsilon (X_n)|\le  |\mathsf E f(X_n)-\mathsf E f^\varepsilon (X_n)|+
|\mathsf E f^\varepsilon (X_n)-\mathsf E g^\varepsilon (X_n)|\\
& \le \varepsilon+C\mathsf P(\varepsilon^{1/2} |X_n|\ge 1)\le \varepsilon+C\varepsilon \mathsf E X_n^2\le\varepsilon+
C\varepsilon\sum_{j=0}^{n-1}\frac{\overline\lambda_j^2\sigma_{j+1}^2}{s_n^2}=(1+C\Lambda^2)\varepsilon.
\end{align*}
From this estimate we obtain the inequalities of the form (\ref{eq:2.6}) with
\[ \mathscr L^\varepsilon=\lim_{n\to\infty}\sup_{\lambda_0^{n-1}\in\mathfrak A^{n-1}}\mathsf E g^\varepsilon\left(X_n\right).\]
Just mentioned inequalities imply that
\begin{equation} \label{eq:2.8}
 \mathscr L:= \lim_{n\to\infty}\sup_{\lambda_0^{n-1}\in\mathfrak A^{n-1}}\mathsf E f(X_n)=\lim_{\varepsilon\to 0}\mathscr L^\varepsilon.
\end{equation}

Denote by $V^\varepsilon$, the viscosity solution of (\ref{eq:1.6}), (\ref{eq:1.7}), corresponding to the terminal condition $g^\varepsilon$ instead of $f$. Since $g^\varepsilon\in C^{2+\alpha}(\mathbb R)$, we have
\begin{equation} \label{eq:2.9}
\mathscr L^\varepsilon=V^\varepsilon(0,0)
\end{equation}
by the result, already proved.

Finally, note, that the convergence  $g^\varepsilon(x)=\chi(\varepsilon^{1/2} x)f^\varepsilon(x)\to f(x)$, $\varepsilon\to 0$ is uniform on compact sets. It follows that
\[ \liminf_{\substack{y\to x\\ \varepsilon\to 0}} g^\varepsilon(y)=\limsup_{\substack{y\to x\\ \varepsilon\to 0}} g^\varepsilon(y)=f(x). \]
Using this fact, by the method of half-relaxed limits, applied above, it is easy to prove that
\begin{equation} \label{eq:2.10}
\lim_{\varepsilon\to 0} V^\varepsilon(0,0)= v(0,0).
\end{equation}

From (\ref{eq:2.8})--(\ref{eq:2.10}) we conclude that $\mathscr L=v(0,0)$. The proof of Theorem \ref{th:1} is complete.

 \bibliographystyle{plain}
 \bibliography{litCLT}

\begin{thebibliography}{10}

\bibitem{AleHer83}
G.~Alefelf and J.~Herzberger.
\newblock {\em Introduction to interval computations}.
\newblock Academic Press, New York, 1983.

\bibitem{BarCap97}
M.~Bardi and I.~Capuzzo-Dolcetta.
\newblock {\em Optimal control and viscosity solutions of
  {H}amilton-{J}acobi-{B}ellman equations}.
\newblock Birkhauser, Boston, 1997.

\bibitem{Bau96}
H.~Bauer.
\newblock {\em Probability theory}.
\newblock Walter de Gruyter, Berlin, 1996.

\bibitem{Boos00}
J.~Boos.
\newblock {\em Classical and modern methods in summability}.
\newblock Oxford University Press, New York, 2000.

\bibitem{CafSte08}
L.A. Caffarelli and U.~Stefanelli.
\newblock A counterexample to ${C}^{2,1}$ regularity for parabolic fully
  nonlinear equations.
\newblock {\em Commun. Part. Diff. Eq.}, 33(7):1216--1234, 2008.

\bibitem{CheGigGot91}
Y.-G. Chen, Y.~Giga, and S.~Goto.
\newblock Remarks on viscosity solutions for evolution equations.
\newblock {\em Proc. Japan Acad., Ser. A}, 67(10):323--328, 1991.

\bibitem{DieFriObe14}
J.~Diehl, P.K. Friz, and H.~Oberhauser.
\newblock Regularity theory for rough partial differential equations and
  parabolic comparison revisited.
\newblock In {\em Stochastic Analysis and Applications 2014}, pages 203--238.
  Springer, 2014.

\bibitem{HuZhou15}
Z.C. Hu and L.~Zhou.
\newblock Multi-dimensional central limit theorems and laws of large numbers
  under sublinear expectations.
\newblock {\em Acta Math. Sin.}, 31(2):305--318, 2015.

\bibitem{KamPelVaz91}
S.~Kamin, L.A. Peletier, and J.L. Vazquez.
\newblock On the {B}arenblatt equation of elasto-plastic filtration.
\newblock {\em Indiana U. Math. J.}, 40(4):1333--1362, 1991.

\bibitem{Kry82}
N.V. Krylov.
\newblock Boundedly nonhomogeneous elliptic and parabolic equations.
\newblock {\em Math. USSR Izv.}, 20(3):459--492, 1983.

\bibitem{LiShi10}
M.~Li and Y.F. Shi.
\newblock A general central limit theorem under sublinear expectations.
\newblock {\em Sci. China Math.}, 53(8):1989--1994, 2010.

\bibitem{Peng07b}
S.~Peng.
\newblock {$G$}-expectation, {$G$}-{B}rownian motion and related stochastic
  calculus of {I}t{\^o} type.
\newblock In {\em Stochastic analysis and applications: The Abel Symposium
  2005}, pages 541--567. Springer, 2007.

\bibitem{Peng07}
S.~Peng.
\newblock Law of large numbers and central limit theorem under nonlinear
  expectations.
\newblock {P}reprint arXiv:math/0702358 [math.PR], 8 pages, 2007.

\bibitem{Peng08}
S.~Peng.
\newblock A new central limit theorem under sublinear expectations.
\newblock {P}reprint arXiv:0803.2656 [math.PR], 25 pages, 2008.

\bibitem{Peng08c}
S.~Peng.
\newblock Nonlinear expectations and stochastic calculus under uncertainty.
\newblock {P}reprint arXiv:1002.4546 [math.PR], 149 pages, 2010.

\bibitem{Pet66}
G.M. Petersen.
\newblock {\em Regular matrix transformations}.
\newblock McGraw-Hill, London, 1966.

\bibitem{Pur67}
L.E. Pursell.
\newblock Uniform approximation of real continuous functions on the real line
  by infinitely differentiable functions.
\newblock {\em Math. Mag.}, 40(5):263--265, 1967.

\bibitem{Rok15}
D.B. Rokhlin.
\newblock Central limit theorem under uncertain linear transformations.
\newblock {P}reprint arXiv:1505.01084v2 [math.PR], 11 pages, 2015.

\bibitem{YongZhou99}
J.~Yong and X.Y. Zhou.
\newblock {\em Stochastic controls: {H}amiltonian systems and {H}{J}{B}
  equations}.
\newblock Springer, New York, 1999.

\bibitem{ZhaChe14}
D.~Zhang and Z.~Chen.
\newblock A weighted central limit theorem under sublinear expectations.
\newblock {\em Commun. Stat. Theory}, 43(3):566--577, 2014.

\end{thebibliography}





\end{document}